\documentclass{amsart}

\usepackage{amsmath,amsthm,amssymb}

\newtheorem{thm}{Theorem}
\newtheorem{lem}[thm]{Lemma}
\newtheorem{pro}[thm]{Proposition}

\def\Diff{\textrm{Diff}}
\def\Homeo{\textrm{Homeo}}
\def\id{\textrm{id}}
\def\R{{\mathbb R}}

\title[Uniform approximation of homeomorphisms]{Uniform approximation of homeomorphisms \\ by diffeomorphisms}
\author{Stefan M\"uller}
\email{stefanm@illinois.edu}
\address{University of Illinois at Urbana-Champaign, Urbana, IL 61801, USA}
\address{Korea Institute for Advanced Study, Seoul 130--722, Republic of Korea}

\subjclass[2010]{57R12, 57Q55, 58C35, 28D05}

\begin{document}
\maketitle

\begin{abstract}
We prove that a compactly supported homeomorphism of a smooth manifold of dimension $n \ge 5$ can be approximated uniformly by compactly supported diffeomorphisms if and only if it is isotopic to a diffeomorphism.
If the given homeomorphism is in addition volume preserving, then it can be approximated uniformly by volume preserving diffeomorphisms.
\end{abstract}

\section{Introduction} \label{sec:intro}
A basic problem in differential topology is to detect if a given homeomorphism of a smooth manifold $M$ can be approximated uniformly by diffeomorphisms.
This question dates back at least as far as \cite{alexander:spt32}.
A solution was found independently by J.~R.~Munkres \cite{munkres:osp60,munkres:hos65} and M.~W.~Hirsch \cite{hirsch:ots63} in the form of an obstruction theory, and by E.~H.~Connell \cite{connell:ash63} for $\R^n$ via an intermediate uniform approximation by PL homeomorphisms.
In particular, a uniform approximation by diffeomorphisms is always possible if the dimension of $M$ is $n \le 3$ \cite{munkres:osp60}, but not necessarily in dimension $4$ \cite{donaldson:agt83}.
In contrast, continuous maps can always be approximated uniformly by smooth maps \cite{hirsch:dt94}.

In this short note, we combine the above results to prove the following theorem.

\begin{thm} \label{thm:approx}
Let $M$ be a closed smooth manifold of dimension $n \ge 5$.
A homeomorphism $\varphi$ of $M$ can be approximated uniformly by diffeomorphisms if and only if $\varphi$ is isotopic to a diffeomorphism.
\end{thm}

See section~\ref{sec:prelim} for details.
We first reduce the theorem to the case of the unit ball in $\R^n$ in section~\ref{sec:local-global}, and then prove the local result in section~\ref{sec:proof-main-lemma}.
Homeomorphisms of non-compact manifolds and manifolds with boundary are discussed in section~\ref{sec:open}.
Not too surprisingly, the problem remains unsolved in dimension $n = 4$.

One can then ask whether a homeomorphism which preserves a given volume form can be approximated uniformly by diffeomorphisms (this problem is mentioned for example in \cite{oxtoby:mph41}), and if so, can it also be approximated by volume preserving diffeomorphisms.
The second part of this question was answered in the affirmative, independently by Y.-G.~Oh \cite{oh:cmp06} and J.-C.~Sikorav \cite{sikorav:avp07}: if a volume preserving homeomorphism can be approximated uniformly by diffeomorphisms, then it can also be approximated uniformly by volume preserving diffeomorphisms.
This is of particular interest in dimension two in the context of $C^0$-symplectic topology \cite{oh:cmp06,mueller:ghh07,mueller:ghh08,mueller:thesis08}, and was the starting point of the author's continued interest in this question.
Interesting variants of this problem have also been studied, e.g.\ approximation by almost everywhere diffeomorphisms \cite{oxtoby:mph41} or in measure (in connection with Lusin's theorem) \cite{white:aom69,white:aom74}.

\section{Preliminaries} \label{sec:prelim}
Let $M$ be a closed smooth manifold of dimension $\dim M \ge 5$.
Without loss of generality, we may assume that $M$ is connected. 
Denote by $\Homeo (M)$ and $\Diff (M)$ the group of homeomorphisms and the subgroup of diffeomorphisms of $M$, respectively.
The group $\Homeo (M)$ is equipped with the uniform (or compact-open) topology.
It is metrizable by the uniform distance of homeomorphisms or by the so-called $C^0$-distance (that is, the uniform distance of homeomorphisms and of their inverses).
Only the latter is complete.
However, if a sequence of homeomorphisms converges uniformly to another homeomorphism (in other words, if a limit exists), then it also converges in the $C^0$-metric.

Recall that any $C^1$-diffeomorphism can be approximated uniformly by $C^\infty$-diffeomorphisms \cite[Theorem 2.7]{hirsch:dt94}. 
Thus for the purposes of this paper, we do not need to distinguish between diffeomorphisms that are of class $C^1$ and $C^\infty$.

An isotopy is a continuous map $\Phi \colon [0,1] \times M \to M$ such that $\varphi_t = \Phi (t, \cdot)$ is a homeomorphism for each $t$; two homeomorphisms $\varphi$ and $\psi$ are said to be isotopic if there exists an isotopy $\Phi$ with $\Phi (0, \cdot) = \varphi$ and $\Phi (1, \cdot) = \psi$.

\section{From local to global} \label{sec:local-global}
The following main lemma is a local version of Theorem~\ref{thm:approx}, and is due to Munkres \cite{munkres:osp60}, Connell \cite{connell:ash63}, and R.~H.~Bing \cite{bing:she63}.

\begin{lem}[Munkres, Connell, Bing] \label{lem:main}
Let $n \ge 5$, and $\varphi \colon B^n \to B^n$ be a homeomorphism of the open unit ball in $\R^n$ that is the identity near the boundary of $B^n$, and is isotopic to the identity through an isotopy that fixes pointwise a neighborhood of the boundary.
Then $\varphi$ can be approximated uniformly by diffeomorphisms of $B^n$ that are the identity near the boundary of $B^n$.
\end{lem}

A proof will be given in the next section.
The only subtlety that is not stated explicitly in the cited references is that the diffeomorphisms that approximate the given homeomorphism can be chosen to be equal to the identity near the boundary, and therefore extend to diffeomorphisms of an ambient manifold.
Aside from a few well-known classical theorems in (differential) topology, this last fact is the main ingredient in our proof of Theorem~\ref{thm:approx}.

Assuming the main lemma, we will give a proof of Theorem~\ref{thm:approx}.

\begin{pro} \label{pro:iso-to-id}
Let $\varphi \in \Homeo (M)$ be a homeomorphism that is isotopic to the identity.
Then $\varphi$ can be approximated uniformly by diffeomorphisms.
\end{pro}

\begin{proof}
Choose a finite cover $\mathcal U$ of $M$ by subsets that are diffeomorphic to open balls in $\R^n$.
By Corollary~1.3 in \cite{edwards:dse71}, $\varphi$ can be fragmented as a finite composition $\varphi = \varphi_m \circ \cdots \circ \varphi_1$, with each $\varphi_i \in \Homeo (M)$ isotopic to the identity by an isotopy that is compactly supported in a member of the cover $\mathcal U$.
By Lemma~\ref{lem:main}, each $\varphi_i$ can be approximated uniformly by diffeomorphisms. 
The proposition follows.
\end{proof}

\begin{proof}[Proof of Theorem~\ref{thm:approx}]
Suppose that $\varphi$ is isotopic to a diffeomorphism $\psi$.
Then the homeomorphism $\psi^{-1} \circ \varphi$ satisfies the hypothesis of Proposition~\ref{pro:iso-to-id}, and therefore can be approximated uniformly by a sequence $\psi_i$ of diffeomorphisms.
But then the sequence $\psi \circ \psi_i$ converges uniformly to $\varphi$.

Conversely, suppose $\varphi$ can be approximated uniformly by diffeomorphisms.
The path components of $\Homeo (M)$ are locally contractible by \cite{cernavskii:lcg69} or \cite{edwards:dse71}, so if $\psi$ is sufficiently close to $\varphi$, then the two homeomorphisms are isotopic.
\end{proof}

\section{Proof of the main lemma} \label{sec:proof-main-lemma}
The main lemma is essentially Theorem~6 in \cite{connell:ash63}, which is stated there with the hypothesis $n \ge 7$.
By the first remark of section~5 of that paper, see also \cite{bing:she63} or \cite[Chapter 4.11]{rushing:te73}, the theorem is also true for $n \ge 5$.
Although this local version of Theorems~3 and 4 in \cite{connell:ash63} is not stated explicitly in \cite[Chapter 4.11]{rushing:te73}, it can be deduced from the previous results of that chapter verbatim as in \cite{connell:ash63}.
The only subtlety not stated explicitly in any of the above references is that with the present hypotheses, the approximating diffeomorphisms can be chosen to equal the identity near the boundary of the ball, and thus extend to diffeomorphisms of the ambient manifold $M$.
The latter however follows from a careful inspection of the proofs of Theorem~6 in \cite{connell:ash63} and Theorem~5.7 in \cite{munkres:osp60}.

\begin{proof}
Extend $\varphi$ to a homeomorphism of $\R^n$ by $\varphi = \id$ outside $B^n$.
Let $\epsilon > 0$.
Then by \cite[Theorem~6]{connell:ash63}, there exists a homeomorphism $\psi$ of $\R^n$ such that $\psi = \id$ on the complement of the open unit ball $B^n$, the restriction of $\psi$ to $B^n$ is a diffeomorphism, and the uniform distance between $\varphi$ and $\psi$ is less than $\epsilon$.
It only remains to show that $\psi$ can be chosen so that it is compactly supported in the interior of the open ball $B^n$, and in particular extends to a diffeomorphism of $\R^n$. 

The final step in the approximation of $\varphi$ by a diffeomorphism is an application of Theorem~5.7 and (the arguments in the proof of) Theorem~6.2 in \cite{munkres:osp60} (and the universal coefficient theorem, cf.\ the proof of Theorem~6.5 in \cite{munkres:osp60}), see the proofs of Theorems~4 and 6 in \cite{connell:ash63}.
If the PL homeomorphisms $f \colon B^n \to B^n$ in \cite[Theorem~6]{connell:ash63} is the identity in a neighborhood of the boundary of the ball, we may delete part of $L$ (the $(n -1)$-skeleton of a smooth triangulation of the open unit ball) in a neighborhood of the boundary without destroying the hypotheses of the theorem.
Then by the last statement of the conclusion, the diffeomorphism that uniformly approximates $f$ may be chosen equal to $f$ and hence the identity in a neighborhood of the boundary of the ball $B^n$.

The intermediate step in the proof is an approximation of $\varphi$ by the PL homeomorphism $f$ above.
This approximation is carried out in sections~3 and 4 of \cite{connell:ash63}.
For the remainder of the present proof, all references to theorems and lemmas are to results in that paper, and for the readers' convenience, the notation is the same as in Connell's paper \cite{connell:ash63} whenever possible.

In the proof of Theorem~6 with $g = \varphi$, one may choose the homeomorphism $h_2 = h_1^{-1}$.
Then $h_2 \circ g \circ h_1 \colon \R^n \to \R^n$ is compactly supported.
The main lemma follows once we show that the PL homeomorphism $h$ in the proof of Theorem~6 can also be chosen to be compactly supported.

This homeomorphism $h$ in turn is a modification of the homeomorphism $h_2 \circ g \circ h_1$ of the preceding paragraph via Theorem~3; more precisely, $h = (h_2 \circ g \circ h_1) \circ f_2 \circ f_1^{-1}$, where $f_1$ and $f_2$ are constructed in the proof of Theorem~3.
(By conjugating by a translation if necessary, the homeomorphism $h_2 \circ g \circ h_1$ itself satisfies the conclusion of the first paragraph in the proof of Theorem~3.)
Therefore it is enough to prove that $f_1$ and $f_2$ can be chosen so that $f_2 \circ f_1^{-1} = \id$ outside a compact subset of $\R^n$.
This can be accomplished by a close inspection of Lemma~4.

Indeed, the PL structures $T$ and $T_1$ agree outside a compact subset of $\R^n$, since the homeomorphism $h_2 \circ g \circ h_1$ is compactly supported.
The two homeomorphisms $f_1$ and $f_2$ are constructed in Lemma~4 as modifications of the same homeomorphism $h$; they are infinite compositions built from repeated applications of Theorem~2.
From the proof of Theorem~2 it follows that each iteration $g_k = f'_k \circ \cdots \circ f'_1$ is obtained from the previous homeomorphism $g_{k -1}$ via composition on the right with another homeomorphism $h'$.
(The homeomorphisms $f'_k$ here are called $f_k$, and the homeomorphism $h'$ is called $h$, in the original proofs.
The notation had to be altered in these two places to avoid an unfortunate overlap in notation among the cited results within \cite{connell:ash63}.)

The above iterated homeomorphisms are piecewise linear with respect to the two PL structures $T$ and $T_1$.
For $k$ sufficiently large, the next step in the induction process modifies the previous step only outside a ball of large radius, where $T$ and $T_1$ coincide.
Thus the homeomorphism $h'$ in the proof of Theorem~2 can be chosen the same in the induction process for $T$ and for $T_1$ in each step for $k$ sufficiently large.
(Again, in the original paper this homeomorphism $h'$ is called $h$.)
As a result, the two iterated homeomorphisms $f_1$ and $f_2$ in the proof of Theorem~3 differ only by a finite number of initial iterations, and all other factors cancel out in the expression $f_2 \circ f_1^{-1}$.

By item (2') on page 333, the composition of the initial iterations is the identity outside a large enough ball.
Therefore $f_2 \circ f_1^{-1}$ is the identity outside a sufficiently large compact subset of $\R^n$.
This completes the proof of the the main lemma.
\end{proof}

\section{Non-compact manifolds and manifolds with boundary} \label{sec:open}
In this section, the manifold $M$ is allowed to be non-compact and to have non-empty boundary.
All homeomorphisms of $M$ are assumed to have compact support in the interior of $M$, and the group $\Homeo (M)$ is equipped with the usual direct limit topology over compact codimension zero submanifolds $K$ of the interior of $M$.
With these modifications, Theorem~\ref{thm:approx} applies to open manifolds:

Let $\varphi$ be a homeomorphism of $M$.
After replacing $M$ by a compact codimension zero submanifold $N$ of its interior such that the interior of $N$ contains the support of $\varphi$, we may assume that $M$ is a compact manifold with non-empty boundary.

Suppose that $\varphi$ is isotopic to a diffeomorphism; the isotopy is supported in a compact subset $K$ of the interior of $M$.
Therefore we may assume without loss of generality that $\varphi$ is isotopic to the identity (cf.\ the proof of Theorem~\ref{thm:approx}).
Then $\varphi$ can be fragmented into homeomorphisms with small support as in the proof of Proposition~\ref{pro:iso-to-id}.
In this case, one needs a relative version of the Edwards-Kirby Corollary~1.3, in which a collar neighborhood of the boundary remains fixed by all deformations; see Section~6 in \cite{edwards:dse71}.
The resulting homeomorphisms $\varphi_i$ are compactly supported in open balls that do not intersect the boundary of $M$.
Thus Lemma~\ref{lem:main} applies, and $\varphi$ can be approximated by a sequence of diffeomorphisms that are compactly supported in a common compact subset of the interior of $M$.

Conversely, suppose that $\psi$ is a diffeomorphism that lies in a small neighborhood of the homeomorphism $\varphi$, and let $N$ be a compact subset of the interior of $M$ that is the closure of an open submanifold of $M$ and contains the support of both $\varphi$ and $\psi$.
Then $\varphi$ and $\psi$ can be viewed as homeomorphisms of $N$ that are the identity on the boundary.
Since the group of homeomorphisms of $N$ that are the identity on the boundary is again locally contractible \cite[Section~7]{edwards:dse71} or \cite{cernavskii:lcg69}, $\varphi$ is isotopic to $\psi$ provided that the two homeomorphism are within $\epsilon$ of one another for sufficiently small $\epsilon > 0$.

For homeomorphism that preserve volume, the only necessary observation is that if $M$ is non-compact, then the volume form in question is assumed to be finite on compact subsets, and thus induces a Radon measure on $M$.
Then the proofs in \cite{oh:cmp06,sikorav:avp07} of uniform approximation by volume preserving diffeomorphisms apply to non-compact manifolds and manifolds with non-empty boundary.
Their result continuous to hold for homeomorphisms and diffeomorphisms that preserve the measure induced by a volume density rather than a genuine volume form on a possibly non-orientable manifold.

\bibliography{uniform}
\bibliographystyle{amsalpha}

\end{document}